\documentclass{notices}
\usepackage{amsfonts,amssymb,amsmath,amscd}

\usepackage{color}

\usepackage{tikz}
	\usetikzlibrary{cd}

\newcommand{\cat}[1]{\ensuremath{\mathsf{#1}}}
\newcommand{\op}{\ensuremath{^\mathrm{op}}}

\DeclareMathOperator{\Ring}{\cat{Ring}}
\newcommand{\Ringop}{\Ring\op}
\DeclareMathOperator{\cRing}{\cat{cRing}}
\newcommand{\cRingop}{\cRing\op}
\DeclareMathOperator{\Set}{\cat{Set}}
\DeclareMathOperator{\Top}{\cat{Top}}

\DeclareMathOperator{\Frac}{Frac}
\DeclareMathOperator{\rad}{rad}

\DeclareMathOperator{\Space}{\cat{Space}}
\DeclareMathOperator{\Alg}{\cat{Alg}}

 \DeclareMathOperator{\pSpec}{\mathit{p}-Spec}

\DeclareMathOperator{\Mod}{\cat{Mod}}
\DeclareMathOperator{\lMod}{\! \mathrm{-} \!\Mod}

\DeclareMathOperator{\Spec}{Spec}
\renewcommand{\O}{\mathcal{O}}
\DeclareMathOperator{\Max}{Max}

\DeclareMathOperator{\Proj}{Proj}
\DeclareMathOperator{\Idpt}{Idem}

\DeclareMathOperator{\Hom}{Hom}
\DeclareMathOperator{\End}{End}

\newcommand{\m}{\mathfrak{m}}
\newcommand{\p}{\mathfrak{p}}
\newcommand{\PP}{\mathcal{P}}

\newcommand{\Z}{\mathbb{Z}}

\newcommand{\C}{\mathbb{C}}
\newcommand{\M}{\mathbb{M}}
\newcommand{\R}{\mathbb{R}}

\numberwithin{equation}{section}
\newtheorem{theorem}[equation]{Theorem}

\title{
A tour of noncommutative spectral theories
}

\author{
  Manuel Reyes
  \affil{
    Manuel Reyes is an associate professor of mathematics at the University of California, Irvine, who can be reached via email at manny.reyes@uci.edu. 
    His work on this article was supported by NSF grant DMS-2201273.
    }
}

\begin{document}

\maketitle

\section{Introduction}


Rings are ubiquitous throughout mathematics. They are home to our earliest numerical manipulations, since all of our familiar number systems form rings. Collections of polynomials, matrices, and continuous functions are just a few common sources of examples that are familiar to all students and practitioners of mathematics. 
Given the wide array of examples, one interesting problem is to seek invariants of rings that can help separate and classify them. For example, there are several numerical invariants of rings and algebras in the form of various dimensions (e.g., Krull, global, Gelfand-Kirillov) that can be assigned to certain rings or algebras.

A particularly effective invariant for commutative rings is the \emph{prime spectrum}: the set of prime ideals of a commutative ring. Originating in number theory as a generalization of prime numbers, this invariant later took on a deeper meaning as it established a strong connection between commutative algebra and geometry. 
Many undergraduate students learn about prime ideals in commutative rings in their introductory courses on rings and fields, so it is also quite widely publicized.
Less well known are the various ways one can extend the notion of a prime ideal to noncommutative rings. 
Inspired by the success of prime ideals in capturing the structure of commutative rings, algebraists have produced a number of different constructions that each restrict to the same prime spectrum for commutative rings.

This paper is a brief survey of various noncommutative generalizations of the spectrum. After a short discussion of sources of noncommutative rings (Subsection~\ref{sub:examples}), we begin with a reminder of the spectrum for commutative rings in Section~\ref{sec:commutative}. We then discuss several ways that the spectrum has been generalized to noncommutative rings in Section~\ref{sec:noncommutative}. 
Finally, Section~\ref{sec:functors} focuses on one of this author's favorite questions: how should we generalize the spectrum to noncommutative rings if we wish for it to define a \emph{functor?} The answers here are only negative so far. We outline a proof of the major obstruction and its surprising connection to the logic of quantum mechanics. The paper concludes with a short discussion of future prospects for noncommutative spectrum functors.

\subsection{Basic examples of noncommutative rings}
\label{sub:examples}

Before proceeding to our discussion of the spectrum of a ring, we offer a high-speed overview of some contexts where noncommutative rings occur. We cannot hope to provide a full introduction to noncommutative ring theory in this brief article, so we highly recommend that interested readers consult the article~\cite{Walton} for a more detailed  survey from a similar perspective or the textbook~\cite{Lam} for a formal introduction to the subject.

Throughout this paper we assume that every ring is associative and has a multiplicative identity. 
We follow the convention that ``noncommutative ring theory'' refers to the study of rings that may or may not be commutative, which is to say that it means ring theory in general.

Historically,  the first noncommutative ring was Hamilton's algebra of \emph{quaternions}, presented with generators and relations as 
\[
\mathbb{H} = \R \langle i,j,k \mid i^2 = j^2 = k^2 = ijk = -1\rangle.
\]
This forms a \emph{division ring}, meaning that every nonzero element has an inverse. In fact, this was Hamilton's key motivation, to construct a number system extending the complex numbers which still possessed a good notion of ``division.'' The great surprise was that, in order to achieve this, he needed to sacrifice commutativity!

Noncommutative algebras occur for natural reasons in linear algebra and related subjects. For any field $k$, the \emph{matrix algebra} $\M_n(k)$ is isomorphic to the ring of linear operators on the vector space $V = k^n$. Noncommutativity of matrices can thus be viewed a consequence of the fact that order of operations typically makes a difference for linear transformations, as it does in the real world. 
Operators on infinite-dimensional vector spaces similarly form noncommutative algebras. In the infinite-dimensional setting, it is often helpful to use tools outside of pure algebra to tame such large algebras. Using tools from functional analysis leads to the study of operator algebras (often acting on a Hilbert space), while tools from topological algebra can give insight into the structure of endomorphisms of infinite-dimensional spaces over arbitrary fields.

Another broad source of noncommutative algebras are those arising from \emph{formal} constructions. These yield rings whose multiplication is not commutative for the simple reason that strings of symbols form ``words'' that are altered when the ``letters'' are rearranged. Chief among these are the free algebras $k \langle x_1, \dots, x_n \rangle$ over a field $k$, whose elements are polynomials with coefficients in $k$ whose terms are formal words in the indeterminates $\{x_i\}$. A closely related construction comes from the theory of quivers (i.e., directed graphs). Each quiver $Q$ defines a path algebra $kQ$, whose elements are formal linear combinations of the paths in $Q$, where the product of two paths is given by composition if one path ends at the beginning vertex of the next, and otherwise is equal to zero.

By far, the most mysterious source of noncommutativity arises from the physical world in the form of \emph{quantum physics.} In order to accurately describe subatomic phenomena, physicists were forced to manipulate algebraic expressions in which some pairs of observable quantities do not commute. 
For instance, the position $\hat{x}$ and momentum $\hat{p}$ of a particle in one dimension are famously postulated to satisfy the canonical commutation relations
\[
\hat{x} \hat{p} - \hat{p} \hat{x} = i \hbar 1,
\]
where $1$ represents an identity operator.
In modern formulations of quantum mechanics, observables are represented by linear operators on a complex Hilbert space. (In many cases, such as the relations above, these operators are required to be unbounded and consequently may only be partially defined operators.)
This commutation relation can be viewed as a \emph{deformation} of the ordinary (trivial) commutation relation between position and momentum of a classical particle, where $\hat{x}\hat{p} - \hat{p} \hat{x} = 0$.
This has inspired a number of other mathematical investigations into deformations of algebras~\cite[Chapter~II]{NCAG}, where a ``classical'' algebra is deformed to a ``quantum'' family of algebras by adjusting parameters in the relations.

\section{Commutative spectral theory}
\label{sec:commutative}

In this section we let $R$ denote a commutative ring. Recall that an ideal $\p$ of $R$ is \emph{prime} if it is a proper ideal (equivalently, $1 \notin \p$) such that, for all $a,b \in R$,
\[
  ab \in \p \implies a \in \p \mbox{ or } b \in \p.
\]
The \emph{spectrum} of a commutative ring $R$ is the set $\Spec R$ of all of its prime ideals.

An important special case of a prime ideal is that of a \emph{maximal} ideal: this is an ideal $\m$ of $R$ that is proper and is maximal among all proper ideals. A typical textbook proof that a maximal ideal is prime uses the fact that an ideal $I$ of $R$ is maximal (respectively, prime) if and only if the quotient ring $R/I$ is a field (respectively, integral domain).
The set of all maximal ideals of $R$ is the \emph{maximal spectrum} of $R$, which we will denote by
\[
  \Max R \subseteq \Spec R.
\]
For instance, the ring of integers $\Z$ has its maximal ideals of the form $(p) = p\Z$ where $p \geq 2$ is a prime number. The only non-maximal prime ideal is $(0)$. 
Similarly, in the ring $k[x]$ of polynomials over a field $k$ in a single indeterminate, the maximal ideals are generated by polynomials that are irreducible over~$k$. If $k$ is algebraically closed (such as $k = \C$), this means that the maximal ideals are of the form $(x - \alpha)$ for $\alpha \in k$. Again, the only non-maximal prime ideal is $(0)$.
The similarity of these two spectra comes from the fact that each of these rings is a prinicpal ideal domain.


Extending the previous example to two variables, we may ask what is the spectrum of the ring of polynomials $k[x,y]$ in two indeterminates. At this point, geometry becomes an indispensable tool to understand the spectrum.
Suppose that $k$ is algebraically closed. In this case, the maximal ideals of $k[x,y]$ correspond to points in the plane $(\alpha, \beta) \in k^2$ as they are of the form $(x-\alpha, y-\beta)$. Again~$(0)$ is a non-maximal prime ideal, but now there are so-called height-1 primes $\p$ such that $(0) \subsetneq \p \subsetneq \m$ for a maximal ideal $\m$. The height-1 primes are of the form $\p = (f)$ for an irreducible polynomial $f(x,y) \in k[x,y]$. 
More generally, the celebrated Nullstellensatz of Hilbert provides a similar geometric description of $\Spec k[x_1, \dots, x_n]$: there is an inclusion-reversing bijection between the set of prime ideals of the polynomial algebra $k[x_1, \dots, x_n]$ and the \emph{irreducible subvarieties} of the affine algebraic variety $V = k^n$. 

The word ``spectrum'' is used not only in commutative algebra but also in linear algebra and operator theory. The connection between these two concepts is more than an accident.
Let $V$ be a finite-dimensional vector space over $\C$. If $T \colon V \to V$ is a linear operator, then its \emph{spectrum} $\sigma(T) \subseteq \C$ is the set of eigenvalues of $T$:
\[
  \sigma(T) = \{\lambda \in \C \mid T - \lambda I \mbox{ is not invertible}\},
\] 
which is finite and nonempty. On the other hand, the operator $T$ generates a complex subalgebra $\C[T]$ of the ring $\End_\C(V)$ of linear endomorphisms of $V$. The elements of $\C[T]$ have the form of complex polynomials evaluated at the operator $T$. In particular, this algebra is commutative and isomorphic to $\C[T] \cong \C[x]/(p(x))$, where $p(x)$ is the minimal polynomial of $T$. If $T$ is diagonalizable (i.e., there is a basis of $V$ consisting of eigenvectors of $T$) then $p(x) = \prod_{\lambda \in \sigma(T)} (x-\lambda)$. An application of the Chinese Remainder Theorem yields
\[
  \C[T] \cong \left. \C[x] \middle/ \left(\prod\nolimits_{\lambda \in \sigma(T)} x - \lambda\right) \cong \C^{\sigma(T)}.\right.
\]
The prime and maximal ideals of this finite direct product of fields all occur as kernels of the projections onto each coordinate, so that
\[
  \Spec \C[T] = \Max \C[T] \cong \sigma(T).
\]
In this way we have a precise correspondence between the ring-theoretic spectrum and the operator-theoretic spectrum.

\subsection{The spectrum as a bridge from algebra to geometry}
\label{sub:geometry}

A remarkable feature of the spectrum, which is not obvious from its initial definition, is the fact that it has deep geometric significance. We might begin to notice this phenomenon in trying to understand the spectrum of polynomial algebras as described above. Yet the story continues, as we now briefly discuss.

First we note that the spectrum of any commutative ring carries a topological structure. For a commutative ring $R$, the \emph{Zariski topology} on $\Spec R$ is the topology whose closed sets are defined to be those subsets of the form
\[
V(I) = \{\p \in \Spec R \mid I \subseteq \p\}
\]
for any ideal $I$ of $R$. This results in an induced topology on $\Max R \subseteq \Spec R$, so that we may view the maximal spectrum as a subspace of the prime spectrum. 

If $f \colon R \to S$ is a homomorphism of commutative rings, then it is an interesting exercise to verify that we have a function
\begin{align*}
f^* \colon \Spec S &\to \Spec R, \\
\p &\mapsto f^{-1}(\p)
\end{align*}
which is continuous with respect to the Zariski topologies. In this way the spectrum becomes a \emph{contravariant} (i.e., arrow-reversing) \emph{functor} from commutative rings to topological spaces, which we denote as a functor out of the opposite category
\[
\Spec \colon \cRingop \to \Top.
\]
(Note that the assignment $f^*$ above does not generally restrict to a function from $\Max S$ to $\Max R$. But in those good cases where it does, we have an induced continuous map between the maximal spectra.)

The deeper geometric meaning of the spectrum arises from a story that can be told across many different flavors of geometry.
Suppose that $X$ is a space. Often there is a space that plays the role of a number line, which typically represents values in a field $k$. Then the set $F(X)$ of maps from $X$ to the line form a commutative ring that is in fact an algebra over~$k$.

In good cases, the algebra $F(X)$ contains rich information about the space $X$. 
This begins at the level of points:
for each $x \in X$ the set 
\[
	\m_x = \{f \in F(X) \mid f(x) = 0\}
\]
of functions that vanish at $x$ will form a maximal ideal of $F(X)$. Furthermore, the assignment
\begin{align*}
X &\to \Max F(X), \\
x &\mapsto \m_x
\end{align*}
tends to form a \emph{bijection} between the points of $X$ and the maximal ideals of its function algebra.

In the nicest possible cases, this connection goes even deeper to the level of categories. The various spaces and appropriate geometric maps between them form a category $\Space$ of spaces. Then the algebra of functions on a space $X$ is the set of morphisms to the line represented by $k$:
\[
F(X) = \Hom_{\Space}(X,k)
\]
This makes it clear that the assignment $X \mapsto F(X)$ forms a contravariant (i.e., arrow-reversing) functor $F = \Hom_{\Space}(-,k)$. 
The algebras $F(X)$ have some natural properties or extra structure coming from the spaces $X$, so that they reside in a certain category $\Alg$ of $k$-algebras. Then the morphisms between function algebras end up being in bijection with the morphisms between their corresponding spaces: for any spaces $X$ and $Y$ we have
\[
\Hom_{\Space}(X,Y) \cong \Hom_{\Alg}(F(Y), F(X)).
\]

In fact, the functor $F$ gives a \emph{dual equivalence} between $\Space$ and $\Alg$, meaning an arrow-reversing equivalence of categories, which can also be written as an equivalence with an opposite category
\[
F \colon \Space\op \overset{\sim}{\longrightarrow} \Alg.
\]
This means that there is a quasi-inverse functor from algebras to spaces, which in each case happens to be the \emph{maximal spectrum!} More precisely, if $A$ is a commutative algebra in $\Alg$ then the space associated to $A$ has underlying set $\Max A$, with a topology that is either the same as or closely related to the relative Zariski topology discussed above.

Some specific instances of these dualities can be given as follows:
\begin{itemize}
\item $\Space$ is the category of compact Hausdorff spaces, and $\Alg$ is the category of commutative C*-algebras.
\item $\Space$ is the category of affine varieties over an algebraically closed field $k$, and $\Alg$ is the category of commutative, reduced, finitely generated $k$-algebras.
\item $\Space$ is the category of Stone spaces (i.e., compact, Hausdorff, totally disconnected spaces), and $\Alg$ is the category of Boolean rings. In this case, the ``number line'' is the field with two elements $k = \mathbb{F}_2$, equipped with the discrete topology.
\end{itemize}

One might naturally wonder whether the prime spectrum functor allows for a similar duality between commutative rings and certain spaces. This is in fact true thanks to the wizardry of  algebraic geometry! 
Without going into great detail, we mention that the functor sending a commutative ring to a spatial object takes the form $R \mapsto (\Spec R, \O_{\Spec R})$ where $\Spec R$ is the prime spectrum equipped with its Zariski topology, and the extra data $\O_{\Spec R}$ is a \emph{sheaf of commutative rings}, which can be imagined as a special assignment of a commutative ring (of ``nice functions'') to each open set of the spectrum. The resulting category $\Space$ in this case is that of \emph{affine schemes}.

\section{Spectral theories for noncommutative rings}
\label{sec:noncommutative}

How should we generalize prime ideals from commutative rings to noncommutative rings? That is to say, how should we define the spectrum of a noncommutative ring? 
Here we will illustrate that this question does not have a single correct answer, which is a common phenomenon in many situations where one wishes to generalize a concept from commutative to noncommutative algebra.

The creative work of algebraists over the past several decades has produced a number of inequivalent ways to generalize the spectrum from commutative rings to noncommutative rings.  Thus, \emph{noncommutative spectral theory} as a whole in fact encompasses a number of different theories. We will briefly discuss a few such theories, in order to illustrate the variety of ways one can understand the spectrum in the noncommutative world.

\subsection{Prime ideals}

The most straightforward way to generalize primes from commutative rings would be to find a suitable definition of prime ideals in a noncommutative ring, in such a way that it specializes to the same set of ideals for commutative rings.
Obviously, we could use the exact same definition of a prime ideal as we did for a commutative ring. Ring theorists define an ideal $\p$ of a noncommutative ring $R$ to be \emph{completely prime} if it is a proper ideal such that, for all $a, b \in R$, if $ab \in \p$ then either $a \in \p$ or $b \in \p$.
The immediate benefit of this definition is that it is familiar, and that a proper ideal $\p$ is completely prime if and only if the ring $R/\p$ has no zero divisors. 
Sadly, the set of these ideals often does not do much to capture the structure of a noncommutative ring. For instance, a ring as simple the set $\M_n(k)$ of square matrices of order $n \geq 2$ over a field $k$ has no completely prime ideals!

The most famous remedy for this problem 
is often summarized~\cite[Chapter~IV]{Lam} by the adage that we can repair the classical definition by ``replacing elements with ideals.'' An ideal $\p$ of a noncommutative ring $R$ is said to be \emph{prime} if it is proper and, for all ideals $I$ and $J$ of $R$, if $\p$ contains the ideal product $IJ$ (the ideal generated by all products $xy$ with $x \in I$ and $y \in J$), then it contains $I$ or $J$:
\begin{equation}\label{eq:prime product}
IJ \subseteq \p \implies I \subseteq \p \mbox{ or } J \subseteq \p.
\end{equation}
We can restate the above condition in terms of elements in the following way: for any $a,b \in R$,
\begin{equation}\label{eq:prime elements}
aRb \subseteq \p \implies a \in \p \mbox{ or } b \in \p.
\end{equation}
One can show that if $R$ happens to be commutative, then the definition above is equivalent to the usual one.

This definition has proved to be quite useful in ring theory and representation theory. 
For instance, every maximal ideal in a ring is prime in the sense above. Because every nonzero ring has a maximal ideal (assuming the Axiom of Choice), it follows that its spectrum of prime ideals is nonempty.
In addition, suppose $V$ is a simple left module (i.e., irreducible representation) over our ring $R$. If $R$ is commutative, then the annihilator of $V$ is a maximal ideal. For a general noncommutative ring $R$, the annihilator of $R$ may not be maximal, but it is always prime.

Another nontrivial application of prime ideals arises in the celebrated work of Goldie~\cite{Goldie:prime}. 
If $R$ is a left noetherian ring, then for every prime ideal $\p$ of $R$, the ring $R/\p$ can be localized at its set of nonzerodivisors, resulting in a ring isomorphic to a matrix ring over a division ring $\M_n(D)$. This is analogous to the commutative case in which $R/\p$ is an integral domain, whose localization produces a field of fractions.


\subsection{Prime left ideals}

In rings that are not commutative, one typically finds left or right ideals that need not form two-sided ideals.
Thus, we might ask whether certain \emph{left} ideals should be considered prime, instead of just two-sided ideals. 
We could reasonably define a left ideal $\p \subsetneq R$ to be \emph{prime} if~\eqref{eq:prime elements} holds for all $a,b \in R$. As before, one can verify that this is equivalent to the condition~\eqref{eq:prime product} where $I$ and $J$ are assumed to be \emph{left} ideals of $R$. A two-sided ideal is prime as a left ideal if and only if it is a prime ideal in the sense defined above.  Prime left ideals in this sense seem to have been introduced at roughly the same time in~\cite{Koh, Michler}. 

This notion is closely connected to that of a \emph{prime module}. A left $R$-module $M \neq 0$ is \emph{prime} if all nonzero submodules of $M$ have the same annihilator ideal in $R$. This is equivalent to the condition that, for all $r \in R$ and $m \in M$,
\[
rRm = 0 \implies rM = 0 \mbox{ or } m = 0.
\]
These modules arise naturally in the study of injective modules over left noetherian rings~\cite[\S V.4]{Gabriel}, as every finitely generated submodule of an indecomposable injective module over such a ring is prime. 
One can then prove~\cite[Proposition~8.1]{Reyes:cohenkaplansky} that a left ideal $\p \subseteq R$ is prime in the sense above if and only if $M = R/\p$ is a prime module.

There are cases where this definition is too ``permissive'' to be useful. For instance, if $R$ is a simple ring (meaning its only ideals are $0 \subsetneq R$), then it turns out that \emph{every} proper left ideal of $R$ is prime! There is a vast array of simple rings whose structure is quite complex, so this notion of prime left ideal does not provide much insight into their structure.

A more restrictive definition of prime left ideal, which in fact generalizes completely prime ideals, was introduced in~\cite{Reyes:prime}. A left ideal $\p \subsetneq R$ is \emph{completely prime} if, for all $a, b \in R$,
\[
\p b \subseteq \p, \  ab \in \p \implies a \in \p \mbox{ or } b \in \p.
\]
An instructive equivalent characterization of this condition can be stated in terms of the endomorphism ring of the left module $R/\p$: the left ideal $\p$ is completely prime if and only every nonzero endomorphism of $R/\p$ is injective. 
We mentioned above that a ring can easily fail to have completely prime ideals, but this is not the case for completely prime \emph{left} ideals. Indeed, every maximal left ideal is completely prime, and every nonzero ring has completely prime left ideals. 
This comparatively new instance of prime left ideals has already found a few unexpected applications. For instance, completely prime left ideals were used in~\cite{AlonParan} in order to formulate a version of the Nullstellensatz for vanishing sets of polynomials in commuting variables with quaternion coefficients. 

\subsection{Primes from division rings}
\label{sub:Cohn}

If $R$ is a commutative ring, then every prime ideal $\p$ of $R$ is the kernel of a homomorphism to a field: namely, the composite morphism $R \to R/\p \to \Frac(R/\p)$ into the field of fractions of the quotient integral domain, or equivalently, the morphism $R \to R_\p \to R_\p/\p R_\p$ to the residue field of the localization. 
By analogy, one might expect primes of a noncommutative ring $R$ to arise from homomorphisms into division rings. Such a theory was developed by Cohn in~\cite[Chapter~7]{Cohn:free}. 
Cohn's insight was that every homomorphism $R \to D$ where $D$ is a division ring has the effect of inverting not only elements of $R$, but also many \emph{matrices} over $R$. Consequently, it is reasonable to consider localizations of rings as homomorphisms that universally invert sets of matrices over $R$. If $\phi \colon R \to D$ is a ring homomorphism, then the set of all square matrices over $R$ that map to non-invertible matrices over $D$ (applying $\phi$ entrywise) is called the \emph{singular kernel} of $\phi$.  
(If $R$ happens to be commutative, then the singular kernel is neatly described as the set of all square matrices whose determinant lies in the prime ideal $\p = \ker \phi$.)


Given a ring $R$, Cohn provided a characterization of all possible singular kernels of homomorphisms from $R$ to any division ring in terms of \emph{prime matrix ideals}, which we now describe.
Given two square matrices $A$ and $B$ over $R$, the \emph{diagonal sum} is the block-diagonal matrix $A \oplus B = \left(\begin{smallmatrix}A & 0 \\ 0 & B\end{smallmatrix}\right)$. Now suppose that $A$ and $B$ are both $n \times n$ matrices that are identical in all but the first column, say $A = (a \ c_2 \ \dots\  c_n)$ and $B = (b \ c_2 \ \dots \ c_n)$. The \emph{determinantal sum} is then defined as $A \nabla B = (a + b \ c_2 \ \dots \ c_n)$. Similarly, if $A$ and $B$ have the same size and coincide in all but a single row or single column, then a determinantal sum can be defined with respect to that row or column. 
(One way to understand the significance of these operations is that, over a commutative ring $R$, the determinant converts them into the ordinary product and sum: $\det(A \oplus B) = \det(A) \det(B)$ and $\det(A \nabla B) = \det A + \det B$ for appropriately chosen matrices $A$ and $B$.)
Also, an $n \times n$ matrix $A$ is said to be \emph{non-full} if there exists a factorization of the form $A = PQ$ where $P$ is an $n \times r$ matrix and $Q$ is an $r \times n$ matrix with $r < n$. (If $R$ is commutative, then such matrices have zero determinant.)
A set $\mathcal{I}$ of square matrices over $R$ is defined to be a \emph{matrix ideal} if it satisfies the following conditions:
\begin{enumerate}
\item $\mathcal{I}$ contains all non-full matrices over $R$;
\item If $A,B \in \mathcal{I}$ and some determinantal sum $C$ of $A$ and $B$ exists, then $C \in \mathcal{I}$;
\item If $A \in \mathcal{I}$ and $B$ is any square matrix over $R$, then $A \oplus B \in \mathcal{I}$;
\item $A \oplus (1) \in \mathcal{I}$ implies that $A \in \mathcal{I}$.
\end{enumerate}
Finally, a matrix ideal $\mathcal{P}$ is \emph{prime} if it does not contain some (equivalently, any) identity matrix and, for any two square matrices $A$ and $B$ over $R$,
\[
A \oplus B \in \mathcal{P} \implies A \in \mathcal{P} \mbox{ or } B \in \mathcal{P}.
\]

Cohn showed that the singular kernel of any homomorphism from a ring $R$ to a division ring $D$ must be a prime matrix ideal.
Conversely, if $\PP$ is any prime matrix ideal over $R$, there is a \emph{localization} $R_\PP$ which is a ring with a homomorphism $R \to R_\PP$ that universally maps matrices outside of $\PP$ to invertible matrices over $R_\PP$. 
Cohn's localization enjoys the following properties:
\begin{itemize}
\item $R_\PP$ is a local ring (i.e., its quotient modulo the Jacobson radical is a division ring);
\item $\PP$ is the singular kernel of the composite ring epimorphism $R \to R_\PP \to R_\PP/\rad R_\PP$.
\end{itemize}
Taken together, the properties of this localization theory ensure that every homomorphism from $R$ to a division ring factors uniquely through the localization $R_\PP$ at a prime matrix ideal $\PP$.


\subsection{Primes from module categories}
\label{sub:abelian}

One other approach to the spectrum originates in noncommutative algebraic geometry, where we imagine that noncommutative rings $R$ could be viewed as ``coordinate rings'' for elusive \emph{noncommutative spaces.} This obviously draws inspiration from commutative spectral theory as outlined in Section~\ref{sub:geometry}.
In many approaches to noncommutative algebraic geometry, one replaces a space with an abelian category of \emph{sheaves of modules} over a structure sheaf on the space. From this perspective, it is reasonable to ask that the spectrum of a ring $R$ be defined in terms of its module category $R \lMod$. 

One of the earliest indications that the spectrum of a ring is reflected in the module category is due to Matlis~\cite{Matlis}. If $R$ is a commutative ring and $\p \in \Spec R$, then the injective hull $E(R/\p)$ of the quotient module is indecomposable. Matlis proved that if $R$ is noetherian, then every indecomposable injective is isomorphic to one of these injective hulls, so that there is a bijection between $\Spec R$ and the isomorphism classes of indecomposable injective modules in $R \lMod$. 

Thus, if $R$ is a left noetherian ring that is not necessarily commutative, one could reasonably view the isomorphism classes of indecomposable injective objects in $R \lMod$ as a kind of ``left spectrum'' of $R$. This led Gabriel~\cite[p.~383]{Gabriel} to define the \emph{spectrum} of a (locally noetherian) abelian category as the collection of isomorphism classes of indecomposable injective objects. He extended Matlis's theorem to the realm of algebraic geometry by showing that the spectrum of the category of quasi-coherent sheaves on a noetherian scheme $X$ is in bijection with the points of $X$, and that in fact one can fully reconstruct the scheme $X$ (including its structure sheaf) from its category of quasi-coherent sheaves.
 This was later extended by Rosenberg~\cite{Rosenberg:reconstruction} to schemes that are not necessarily noetherian, with the help of a more general notion of spectrum for abelian categories; we also refer readers to the careful treatment in~\cite{Brandenburg:reconstruction}. 

Since then, a number of other possible spectra of abelian categories have been proposed. These are defined in a variety of ways, from equivalence classes of objects in the category to special subcategories (e.g., Serre or localizing subcategories). A number of these have been surveyed in~\cite[Chapter~VI]{Rosenberg:book}. More recent developments include the \emph{atom spectrum} of Kanda~\cite{Kanda}, which is in bijection with the injective spectrum for locally noetherian Grothendieck categories but generalizes well in the absence of the noetherian hypothesis. 

%
%

%

%

\section{Noncommutative spectrum functors}
\label{sec:functors}

We see now that there is no shortage of interesting ways to extend the prime spectrum from commutative rings to noncommutative rings. 
In the case of commutative rings, we have also seen that functoriality of the spectrum is a fundamental property. 
Thus it seems natural to expect that a spectrum construction for noncommutative rings should form a functor. However, we will soon learn that the reality is much more complicated.

To allow for precise discussion, we refer to the following diagram:
\begin{equation}\label{eq:spec diagram}
\begin{tikzcd}
\cRing\op \ar[r, "\Spec"] \ar[d, phantom, sloped, "\subseteq"] & \Top \ar[r, "U"] & \Set \\
\Ring\op \ar[ur, dashrightarrow, "\exists F?"] \ar[urr, dashrightarrow] & &
\end{tikzcd}
\end{equation}
The symbol ``$\subseteq$'' denotes the inclusion of a full subcategory, and $U$ is the forgetful functor to the category of sets. 
The most basic question to ask is whether the functor $\Spec$ on the category of commutative rings has an extension $F$ to the category of all rings. It is of course natural to ask for such a functor to the target category of topological spaces. But in order to construct a topological space, we first need an underlying set, and this problem is already interesting enough to keep us occupied. For this reason we will focus on functors
\[
F \colon \Ringop \to \Set
\]
whose restriction to commutative rings 
\[
F|_{\cRingop} \cong \Spec
\] 
is naturally isomorphic to the Zariski spectrum, which is to say that we have bijections for each commutative ring $C$
\[
F(C) \cong \Spec(C)
\]
that are natural in $C$. 

Are there \emph{any} functors $F$ that give a noncommutative extension of $\Spec$ in the sense above? In fact, there are several options. The most obvious choice would be to take $F(R)$ to be the set of all completely prime ideals of $R$, which forms a contravariant functor as the preimage of a completely prime ideal under a ring homomorphism is again completely prime. One could also define a functor $F$ by taking the composite
\[
\begin{tikzcd}
\Ringop \ar[r] & \cRingop \ar[r, "\Spec"] & \Set \\
R \ar[r, mapsto] & R/[R,R] \ar[r, mapsto] & \Spec R/[R,R]
\end{tikzcd}
\]
where $[R,R]$ denotes the ideal generated by all commutators in $R$. If $C$ is a commutative ring, then $C/[C,C] \cong C$ so that we obtain natural isomorphisms $F(C) \cong \Spec C$.  
For a third example that is not so obvious, we could let $F(R)$ be the set of prime matrix ideals of $R$ as discussed in Section~\ref{sub:Cohn}. One can show that the preimage of a prime matrix ideal under a ring homomorphism is again a prime matrix ideal; intuitively, this is because each prime matrix ideal of a ring $S$ is the singular kernel of a ring homomorphism $S \to D$ for a division ring, which then induces a singular kernel for the composite $R \to S \to D$ along any homomorphism $f \colon R \to S$. In this way, prime matrix ideals form a functor $F \colon \Ringop \to \Set$ in such a way that we have natural bijections $F(R) \cong \Spec R$ for commutative rings $R$.

However, there is a common deficiency for each of these functors. For any field $k$, if we set $R = \M_n(k)$ for $n \geq 2$, then it is a brief exercise to show that $R$ has no completely prime ideals, satisfies $R/[R,R] = 0$, and has no homomorphisms to any division ring. Thus for each of the functors above we have $F(\M_n(k)) = \varnothing$. Surely this is an undesirable situation. 
If we wish for our spectrum functor $F$ to properly reflect the internal structure of a ring, then it should return nontrivial information for such a simple ring as a matrix algebra. In fact, by comparison with the commutative spectrum we would naturally expect that $F(R) \neq \varnothing$ for every nonzero ring $R$.

On first reflection it seems that this is likely the fault of the particular functors $F$ that we have chosen above. Surely we can be more creative in our choice of functor to avoid this problem, can't we? 
For instance, several of the noncommutative spectra described in Section~\ref{sec:noncommutative} assign nonempty sets to every nonzero ring; certainly one of them can be made into a functor? 
If one carefully studies how these various flavors of primes should pull back along a ring homomorphism $f \colon R \to S$, it becomes apparent that there are serious difficulties. (For each of the various spectra, this difficulty with pulling back primes already occurs for the injective homomorphism $k \times k \hookrightarrow \M_2(k)$ given by the diagonal embedding.) 

Much to the shock of this author, it turns out that \emph{every} possible choice of spectrum functor shares this deficiency! The following states the theorem precisely.

\begin{theorem}[\cite{Reyes:obstructing}] \label{thm:empty spec}
Let $F \colon \Ringop \to \Set$ be a functor such that there are natural bijections $F(C) \cong \Spec(C)$ for every commutative ring $C$. Then $F(\M_n(\C)) = 0$ for every $n \geq 3$.
\end{theorem}

In the next two subsections we outline the ideas involved in the proof of this theorem.

\subsection{A universal functor}
How can we prove that such an obstruction holds for \emph{all} possible spectrum functors $F$?
The key observation that makes the proof possible is the fact that among all such functors, there is a \emph{universal} one. Understanding this functor requires an unusual but simple definition. Let $R$ be a ring and let $\p \subseteq R$ be a subset. We say that $\p$ is a \emph{prime partial ideal} if, for every commutative subring $C \subseteq R$, the intersection $C \cap \p$ is a prime ideal of $C$. We then define the \emph{partial spectrum} of $R$ to be the set $\pSpec R$ of all prime partial ideals of $R$. It is a simple exercise to verify that the following hold:
\begin{itemize}
\item If $R$ is commutative, then $\pSpec R = \Spec R$.
\item If $f \colon R \to S$ is a ring homomorphism and $\p \in \pSpec S$, then $f^{-1}(\p) \in \pSpec R$. 
\end{itemize}
Together, this shows that we have a functor
\[
\pSpec \colon \Ringop \to \Set
\]
whose restriction to commutative rings coincides with the usual spectrum functor. 

In what way is this functor universal? Suppose that $F$ is any functor as in the statement of Theorem~\ref{thm:empty spec}, along with a fixed natural isomorphism $\alpha \colon F|_{\cRingop} \overset{\sim}{\longrightarrow} \Spec$. Then there exists a unique natural transformation $\eta \colon F \to \pSpec$ 
whose restriction to $\cRingop$ is equal to $\alpha$. 

This universal property is significant for the following reason. Suppose that we are able to prove $\pSpec(R) = \varnothing$ for some ring $R$. Then, for any functor $F$ as above, by the universal property there exists a function
\[
\begin{tikzcd}
F(R) \ar[r, "\eta_R"] & \pSpec(R) = \varnothing.
\end{tikzcd}
\]
But the only set with a function to the empty set is the empty set itself. It would then follow that $F(R) = \varnothing$ for \emph{every} functor $F$ as above!

\subsection{From partial ideals to quantum colorings}
Thus we have a new strategy: to prove Theorem~\ref{thm:empty spec}, it suffices to show that $\pSpec \M_n(\C) = \varnothing$ for $n \geq 3$. For simplicity we focus on the case $n = 3$. (In fact, the general case reduces to this critical value.) 
The next key observation is that if the partial spectrum of $\M_3(\C)$ were nonempty, then each prime partial ideal would induce a special type of coloring on projection matrices. Indeed, suppose that there exists a prime partial ideal $\p$ of $\M_3(\C)$. 
%
We can define a $\{0,1\}$-valued coloring on the set $\Proj(\M_3(\C))$ of all projection matrices by setting
\begin{equation}\label{eq:coloring}
c_{\p}(P) = \begin{cases}
1, & P \notin \p, \\
0, & P \in \p.
\end{cases}
\end{equation}
This function $c_{\p} \colon \Proj(\M_3(\C)) \to \{0,1\}$ has the following special property: if $P,Q \in \Proj(\M_3(\C))$ are projections that commute under matrix multiplication ($PQ = QP$), then $c_\p(P+Q) = c_\p(P) + c_\p(Q)$ and $c_\p(PQ) = c_\p(P)c_\p(Q)$. 

Note that if $C \subseteq \M_3(\C)$ is a commutative subalgebra, then $\Proj(C)$ forms a Boolean algebra under the operations $P \wedge Q = PQ$, $\neg P = I-P$, and $P \vee Q = \neg(\neg P \wedge \neg Q) = P+Q-PQ$. It follows that this function $c_\p$ restricts to a homomorphism of Boolean algebras $\Proj C \to \{0,1\}$.

At this point in the story, a wonderful coincidence occurs: The functions above were studied decades ago in order to understand the logical foundations of quantum mechanics! Kochen and Specker~\cite{KochenSpecker} investigated a certain kind of hidden variable theory, which is an attempt to reduce our understanding of quantum mechanics to underlying classical variables. 
Viewing $\M_3(\C)$ as the algebra of observables on a quantum system, each projection represents an observable whose possible values are either~0 or~1 (given by the eigenvalues of the operator); in essence, these are ``yes-no questions'' that we can ask about the system. 
A suitably general version of Heisenberg's uncertainty relations~\cite[Proposition~1.4]{Takhtajan} implies that we cannot precisely measure the value of two such observables unless they happen to commute; for this reason, commuting observables are referred to as \emph{commeasurable}.

In this context, the coloring $c = c_\p$ assigns a definite value to all $\{0,1\}$-valued observables in such a way that it restricts to a Boolean algebra homomorphism on subsets of commeasurable observables. Kochen and Specker referred to this as a \emph{morphism of partial Boolean algebras}.
They then proved~\cite[Theorem~1]{KochenSpecker} that no such coloring exists by constructing a finite set of vectors whose corresponding rank-1 projections cannot be colored.

\begin{theorem}[\cite{KochenSpecker}]
There is no coloring $c \colon \Proj(\M_3(\C)) \to \{0,1\}$ that satisfies the conditions $c(P+Q) = c(P) + c(Q)$ and $c(PQ) = c(P)c(Q)$ for all pairs of commeasurable projections $PQ = QP$.
\end{theorem}

The existence of the function $c_\p$ defined above evidently contradicts the Kochen-Specker Theorem, which implies that $\M_3(\C)$ cannot have any prime partial ideals. Thus $\pSpec(\M_3(\C)) = \varnothing$, from which Theorem~\ref{thm:empty spec} follows.

\subsection{Further remarks}

There are a few other natural questions to ask regarding Theorem~\ref{thm:empty spec}. First, what happens in the case of $2 \times 2$ matrices where $n = 2$? In fact, if we choose the functor $F = \pSpec$, it turns out that $F(\M_2(\C))$ is a very large set, of cardinality $2^{2^{\aleph_0}}$! The conclusion of the theorem fails spectacularly in this case. 

A second natural question is whether the theorem can be generalized beyond complex matrices. Is it still true that for other fields or even other rings $R$, every functor $F$ as in the statement of the theorem satisfies $F(\M_n(R)) = \varnothing$ if $n \geq 3$? The answer in this case is affirmative, as shown in~\cite{BMR:KS}. This problem reduces to the universal case $R = \Z$. The proof for the ring of integers follows from an analogue of the Kochen-Specker theorem for colorings of idempotent matrices in $\M_3(\Z)$. We let $\Idpt R$ denote the set of idempotents of any ring $R$.

\begin{theorem}[\cite{BMR:KS}]
There is no coloring $c \colon \Idpt(\M_3(\Z)) \to \{0,1\}$ such that $c(P+Q) = c(P)+c(Q)$ and $c(PQ) = c(P)c(Q)$ for all pairs of commuting idempotents $P$ and $Q$.
\end{theorem}

If a prime partial ideal $\p \in \pSpec(\M_3(\Z))$ exists, then we can again define a coloring $c_\p \colon \Idpt(\M_3(\Z)) \to \{0,1\}$ as in~\eqref{eq:coloring}, contradicting the theorem above. Thus $\pSpec(\M_3(\Z)) = \varnothing$. For any ring $R$, we have a unique homomorphism $\Z \to R$ (given by $m \mapsto m \cdot 1_R \in R$), which extends entrywise to a homomorphism $f \colon \M_3(\Z) \to \M_3(R)$. Now if $F$ is a functor as in Theorem~\ref{thm:empty spec}, we obtain a function
\[
\begin{tikzcd}
F(\M_3(R)) \ar[r, "F(f)"] & F(\M_3(\Z)) \ar[r, "\eta"] &[-.5em] \underbrace{\pSpec(\M_3(\Z))}_{= \varnothing}.
\end{tikzcd}
\]
It follows that $F(\M_3(R)) = \varnothing$ as well!

\subsection{Outlook}


What does Theorem~\ref{thm:empty spec} mean for the future of noncommutative spectral theory? 
On the one hand, it is clear from the variety of ideas surveyed in Section~\ref{sec:noncommutative} that noncommutative spectral theory is alive and well. Indeed, a spectrum does not need to be functorial in order to provide us with interesting information! 
These spectra join many other important invariants that are not functors, such as the center of a ring or the Hochschild cohomology of an algebra over a field. Functoriality is a wonderful property when it holds, but it is not necessary for a construction to be useful.

On the other hand, the functorial correspondences in subsection~\ref{sub:geometry} still provide compelling inspiration for \emph{noncommutative geometry}~\cite{NCAG}, a multifaceted subject which has been a steady source of rich mathematical discovery for several decades now. For those of us who are not willing to give up on the dream of ``noncommutative spaces'' that are dual to categories of noncommutative rings and algebras, what is the path forward?

The author has written about such prospects elsewhere~\cite[Section~1]{Reyes:coalgebra}; only a hint of this is recalled below in the hopes that it will capture the imagination of future contributors to this problem. One interpretation of Theorem~\ref{thm:empty spec} and related no-go theorems is that spaces built out of point sets are ``too commutative'' (i.e., classical) to properly model noncommutative behavior, and in order to build a truly noncommutative algebra-geometry correspondence we need a suitable noncommutative (perhaps even \emph{quantum}) generalization of sets. 
This would amount to finding a category $\mathfrak{S}$ of generalized sets that should ideally have a fully faithful embedding $\Set \hookrightarrow \mathfrak{S}$ from the category of classical sets. To repair the obstruction to~\eqref{eq:spec diagram}, a noncommutative spectrum functor should be a functor $\Sigma \colon \Ringop \to \mathfrak{S}$ such that the diagram
\begin{equation}\label{eq:nc spec}
\begin{tikzcd}
\cRingop \ar[r, "\Spec"] \ar[d, hookrightarrow] & \Top \ar[r, "U"] & \Set \ar[d, hookrightarrow] \\
\Ringop \ar[rr, "\Sigma"] & & \mathfrak{S}
\end{tikzcd}
\end{equation}
commutes. The effectiveness of any such candidate functor $\Sigma \colon \Ringop \to \mathfrak{S}$ could be measured in a variety of ways. Does it assign a nontrivial object to every ring? Could it allow us to recover any of the various noncommutative spectra surveyed in Section~\ref{sec:noncommutative}? For a ring $R$, can $\Sigma(R)$  be enhanced to a complete invariant of $R$, similar to the way that the spectrum of a commutative ring can be equipped with a structure sheaf? Is the resulting category $\mathfrak{S}$ rich enough to allow for the construction of various kinds of noncommutative spaces? 

It is not yet clear how exotic the objects of $\mathfrak{S}$ should be in order to provide a nontrivial spectrum for every noncommutative ring. As discussed in~\cite{Reyes:coalgebra}, some hints point toward such objects being roughly related to \emph{coalgebras}, although in full generality we cannot expect these objects to depend on a single base field. Thus it remains a wide open problem to realize the picture presented in~\eqref{eq:nc spec}. 
In light of the history of creative reinterpretations of the spectrum, we remain hopeful that sooner or later a suitable candidate will be found. We hope that this all-too-brief survey will spark the interest of a few readers who might help to make this a reality. 

\section*{Acknowledgments}

It is my pleasure to thank T.\,Y.~Lam, So Nakamura and Ariel Rosenfield for feedback on a draft of this article.
I am also grateful to the anonymous reviewers and to Han-Bom Moon for several corrections and suggestions to for improvement.

\bibliography{ncspec}

\end{document}